\documentclass[12pt,fleqn]{article}
\usepackage{amsfonts,amssymb,amsmath,amsthm,graphicx,multicol}

\newtheorem{theorem}{Theorem}

\title{Phase Plots of Complex Functions: \\ a Journey in Illustration}

\author{Elias Wegert, TU Freiberg}

\date{July 14, 2010}

\begin{document}

\maketitle

\vspace*{10mm}

\begin{multicols}{2}

\section*{Introduction}

This work was inspired by the recent article ``M\"obius Transformations
Revealed'' by Douglas Arnold and Jonathan Rogness \cite{ArnRog}.
There the authors write:
\smallskip

``{\em Among the most insightful tools that mathematics has developed is the
representation of a function of a real variable by its graph. \ldots\
The situation is quite different for a function of a complex variable.
The graph is then a surface in four dimensional space, and not so easily
drawn. Many texts in complex analysis are without a single depiction of
a function. Nor it is unusual for average students to complete a course in
the subject with little idea of what even simple functions, say trigonometric
functions, `look like'.}''
\smallskip

In the printed literature there are a few laudable exceptions to this rule,
such as the prize-winning ``Visual Complex Analysis'' by Tristan Needham
\cite{Need}, Steven Krantz' textbook \cite{Kra2} with
a chapter on computer packages for studying complex variables, and
the {\sc Maple} based (German) introduction to complex function theory
\cite{ForHof} by Wilhelm Forst and Die\-ter Hoffmann.

But looking behind the curtain, one encounters a different situation which is
evolving very quickly.
Many of us have developed our own techniques for visualizing complex
functions in teaching and research, and one can find many beautiful
illustrations of complex functions on the internet.

This paper is devoted to ``phase plots'', a special tool for visualizing
and exploring analytic functions. Figure~1 shows such a fingerprint
of a function in the complex unit disk.
\begin{center}\label{f.fig1}
\includegraphics[width=0.45\textwidth]{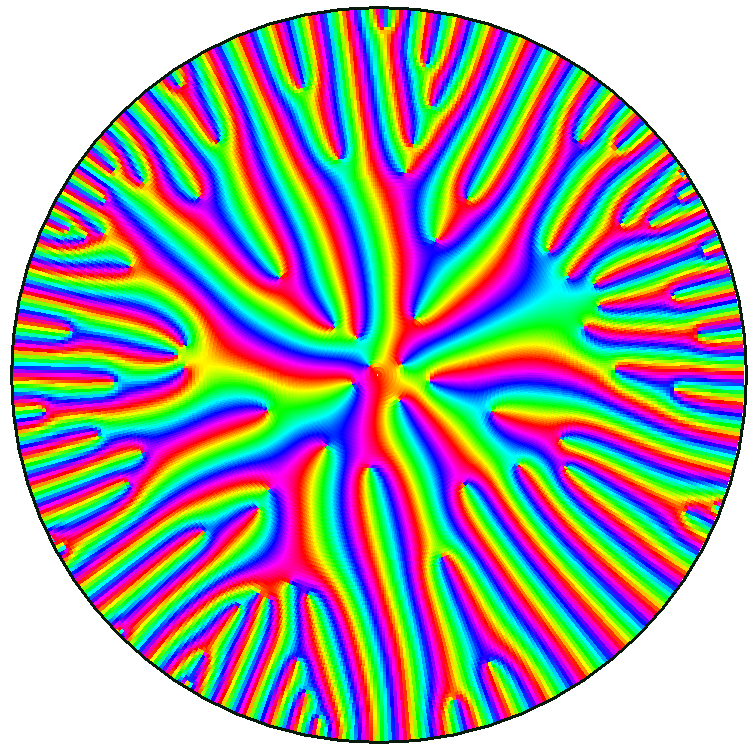} \\
Figure~1: The phase plot of an analytic function in the unit disk
\end{center}
The explanation of this illustration is deferred to a later section
where it is investigated in detail.

Phase plots have been invented independently by a number of people
and it is impossible to give credit to someone for being the first.
Originally, they were mainly used in teaching as simple and effective
methods for visualizing complex functions.
Over the years, and in particular during the process of writing and
rewriting this manuscript, the topic developed its own dynamics and
gradually these innocent illustrations transmuted to sharp tools
for dissecting complex functions.

So the main purpose of this paper is not only to present nice
pictures which allow one to recognize complex functions by their
individual face, but also to develop the mathematical background and
demonstrate the utility and creative uses of phase plots. That they
sometimes also facilitate a new view on known results and may open up new
perspectives is illustrated by a universality property of the
Riemann Zeta function which, in the setting of phase plots, can be
explained to (almost) anyone.

The final section is somewhat special. It resulted from a
self-experiment carried out to demonstrate that phase plots
are sources of inspiration which can help to establish new results.
The main finding is that any meromorphic function is associated
with a {\em dynamical system} which generates a {\em phase flow} on
its domain and converts the phase plot into a {\em phase diagram}.
These diagrams will be useful tools for exploring complex functions,
especially for those who prefer thinking geometrically.

\section*{Visualization of Functions}

The graph of a function $f: D\subset\mathbb{C} \rightarrow
\mathbb{C}$ lives in four real dimensions, and since our imagination is
trained in three dimensional space, most of us have difficulties
in ``seeing'' such an object.%
\footnote{One exception is Thomas Banchoff who visualized four
dimensional graphs of complex functions at \cite{Ban2}.}

Some old books on complex function theory have nice illustrations of
analytic functions. These figures show the {\em analytic landscape}
of a function, which is the graph of its modulus.
\begin{center}\label{f.fig2}
\includegraphics[width=0.4\textwidth]{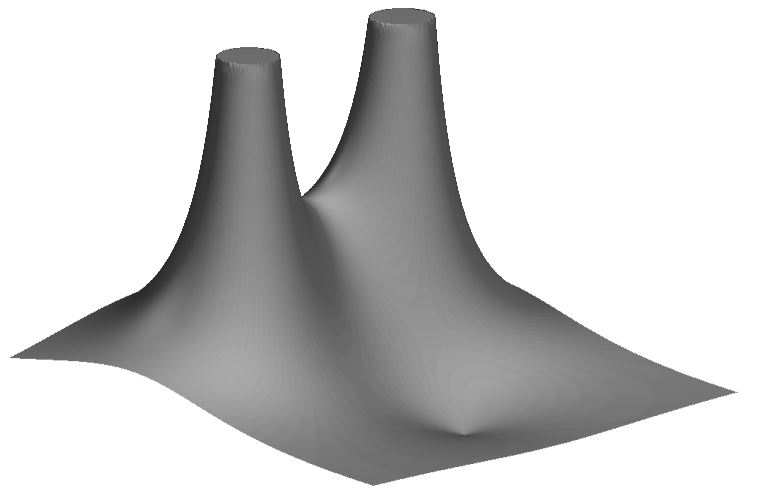} \\
Figure~2: The analytic landscape of $f(z)=(z-1)/(z^2+z+1)$
\end{center}
The concept was not introduced by Johann Jensen in 1912 as sometimes claimed,
but probably earlier by Edmond Maillet \cite{Mai} in 1903 (see also
Otto Reimerdes' paper \cite{Reim} of 1911).

Differential geometric properties of analytic landscapes have been
studied in quite a number of early papers (see Ernst Ullrich
\cite{Ullr} and the references therein). Jensen \cite{Jen} and
others also considered the graph of $|f|^2$, which is a smooth surface.
The second edition of ``The Jahnke-Emde'' \cite{JahEmd} made
analytic landscapes popular in applied mathematics.

Analytic landscapes involve only one part of the function $f$,
its {\em modulus} $|f|$; the {\em argument} $\arg f$ is lost. In the era
of black-and-white illustrations our predecessors sometimes compensated
this shortcoming by complementing the analytic landscape with
lines of constant argument. Today we can achieve this much better using
colors. Since coloring is an essential ingredient of phase plots we
consider it in some detail.

Recall that the argument $\arg z$ of a complex number $z$ is unique up
to an additive multiple of $2\pi$. In order to make the argument
well-defined its values are often restricted to the interval $(-\pi,\pi]$,
or, even worse, to $[0,2\pi)$.
This ambiguity disappears if we replace $\arg z$ with the {\em phase}
$z/|z|$ of $z$.
Though one usually does not distinguish between the notions of ``argument''
and ``phase'', it is essential here to keep these concepts apart.

The phase lives on the complex unit circle $\mathbb{T}$, and points on
a circle can naturally be encoded by {\em colors}. We thus let color serve
as the lacking fourth dimension when representing graphs of complex-valued
functions.
\begin{center}\label{f.fig3}
\includegraphics[width=0.2\textwidth]{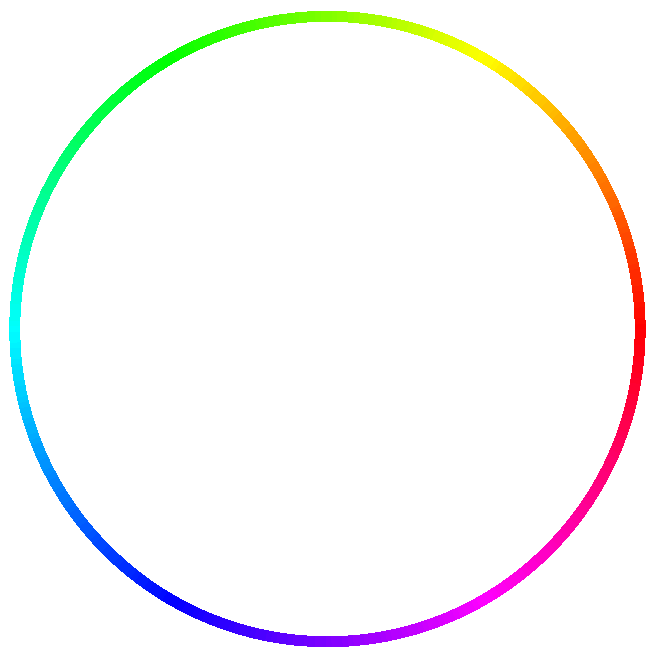} \hfill
\includegraphics[width=0.2\textwidth]{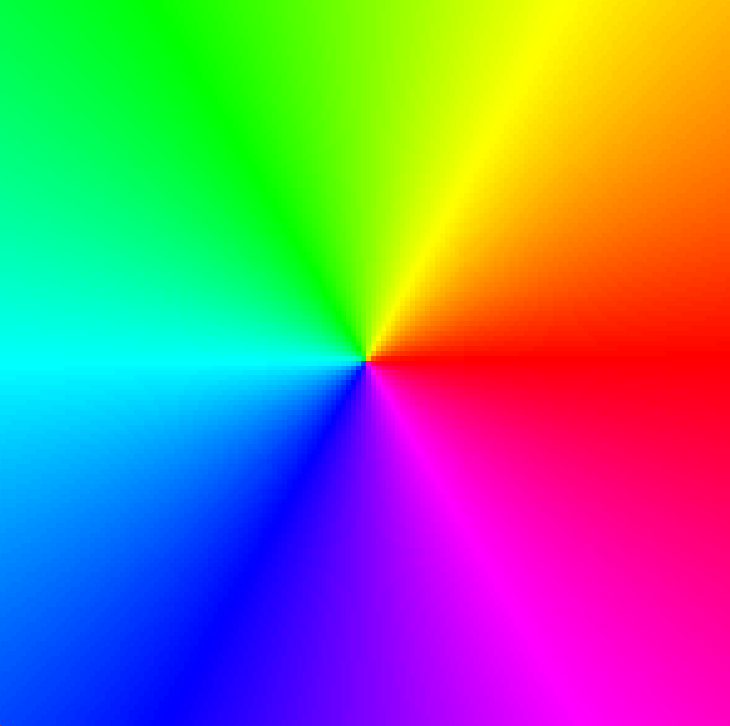} \\
Figure~3: The color circle and the color-coded phase of points close to the
origin
\end{center}
The {\em colored analytic landscape} is the graph of $|f|$ colored according
to the phase of $f$. Since the modulus of analytic functions typically varies
over a wide range one better uses a logarithmic scaling of the vertical axis.
This representation is also more natural since $\log |f|$ and $\arg f$ are
conjugate harmonic functions.
\begin{center}\label{f.fig4}
\includegraphics[width=0.45\textwidth]{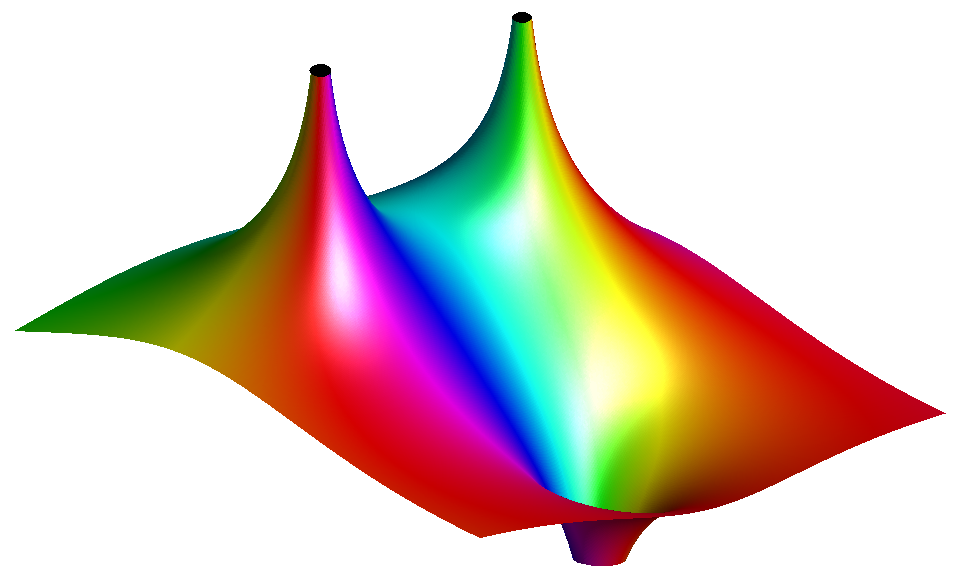}\\
Figure~4: The colored logarithmic analytic landscape of $f(z)=(z-1)/(z^2+z+1)$
\end{center}
Colored analytic landscapes came to life with easy access to computer
graphics and by now quite a number of people have developed software for
their visualization. Andrew Bennett \cite{Benn} has an easy-to-use Java
implementation, and an executable Windows program can be downloaded
from Donald Marshall's web site \cite{Mar}. We further refer to Chapter~12 of
Steven Krantz' book \cite{Kra2}, as well as to the web sites run by
Hans Lundmark \cite{Lund} and Tristan Needham \cite{NeedWeb}.
Very beautiful pictures of (uncolored) analytic landscapes can be found
on the ``The Wolfram Special Function Site'' \cite{Wolf}.%
\footnote{As of June 2009 Wolfram's tool visualizes the analytic landscape
and the argument, but not phase.}

Though color printing is still expensive, colored analytic landscapes
meanwhile also appear in the printed literature (see, for example, the
outstanding mathematics textbook \cite{HoeMath} by Arens et al.\
for engineering students).

With colored analytic landscapes the problem of visualizing complex
functions could be considered solved. However, there is yet
another approach which is not only simpler but even more general.

Instead of drawing a graph, one can depict a function directly on its
domain by color-coding its values, thus converting it to an {\em image}.
Such {\em color graphs} of functions $f$ live in the product of the domain
of $f$ with a color space.

Coloring techniques for visualizing functions have been customary for
many decades, for example in depicting altitudes on maps, but mostly they
represent {\em real valued} functions using a one dimensional color scheme.
It is reported that two dimensional color schemes for visualizing
complex valued functions have been in use for more than twenty years by now
(Larry Crone \cite{Crone}, see Hans Lundmark \cite{Lund}), but they became
popular only with Frank Farris' review \cite{Farr} of Needham's book and its
complement \cite{Far2}. Farris also coined the name ``domain coloring''.

Domain coloring is a natural and universal substitute for the graph of a
function. Moreover it easily extends to functions on Riemann surfaces or
on surfaces embedded $\mathbb{R}^3$ (see Konstantin Poelke and Konrad Polthier
\cite{PoePol}, for instance).

It is worth mentioning that we human beings are somewhat limited with
respect to the available color spaces. Since our visual system has
three different color receptors, we can only recognize colors from
a three--dimensional space.
Mathe\-ma\-ti\-cians of the species {\sf gonodactylus oerstedii}%
\footnote{
This is a species of shrimps which have 12 different photo receptors.}
could use domain coloring techniques to even visualize functions with values
in a twelve--dimensional space (Welsch and Liebmann \cite{WelLie} p.\,268, for
details see Cronin and King \cite{Cro}).

Indeed many people are not aware that natural colors in fact provide us with
an {\em infinite dimensional space} - at least theoretically.
In reality ``color'' always needs a carrier. ``Colored light'' is an
electromagnetic wave which is a mixture of monochromatic components
with different wavelengths and intensities. A simple prismatic piece of
glass reveals how light is composed from its {\em spectral
components}. Readers interested in further information are recommended
to visit the fascinating internet site of Dieter Zawischa \cite{Zaw}.
\begin{center}\label{f.fig5}
\includegraphics[width=0.45\textwidth]{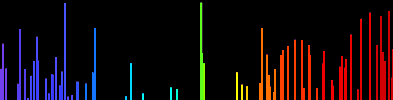} \\
Figure~5: A typical spectrum of polar light
\end{center}
The wavelengths of visible light fill an interval between $375$ nm and $750$
nm approximately, and hence color spectra form an infinite-dimensional space.

How many color dimensions are distinguishable in reality depends on the
resolution of the measuring device. A simple model of the human eye,
which can be traced back to Thomas Young in 1802,
assumes that our color recognition is based on three types of receptors
which are sensitive to red, green, and blue, respectively.

Since, according to this assumption, our visual color space has
dimension three, different spectra of light induce the same visual
impression. Interestingly, a mathematical theory of this effect
was developed as early as in 1853 by Hermann Grassmann,
the ingenious author of the ``Ausdehnungslehre'',
who found three fundamental laws of this so-called {\em metamerism}
\cite{Gra} (see Welsch and Liebmann \cite{WelLie}).

Bearing in mind that the world of real colors is infinite dimensional,
it becomes obvious that its compression to at most three dimensions
cannot lead to completely satisfying results, which explains
the variety of color schemes in use for different purposes.
The two most popular color systems in our computer dominated world
are the RGB (CYM) and HSV schemes.

In contrast to domain colorings which color-code the complete values
$f(z)$ by a two dimensional color scheme, {\em phase plots} display
only $f(z)/|f(z)|$ and thus require just a {\em one dimensional} color space
with a {\em circular topology}.
As will be shown in the next section, they nevertheless contain almost
all relevant information about the depicted analytic or meromorphic
function.

In the figure below the Riemann sphere $\widehat{\mathbb{C}}$
(with the point at infinity on top) is colored using two typical schemes
for phase plots (left) and domain coloring (right), respectively.
\begin{center}\label{f.fig6}
\includegraphics[width=0.20\textwidth]{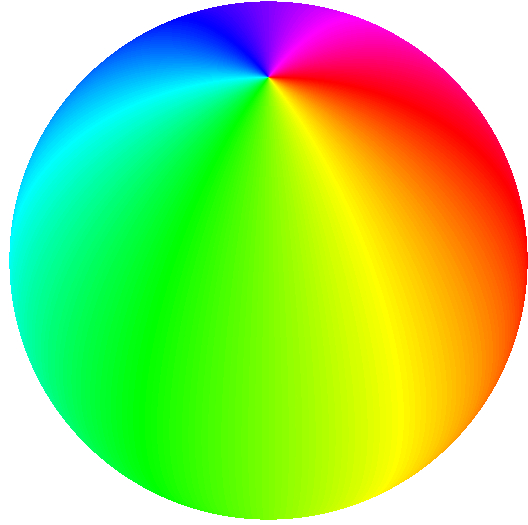}
\hspace{10mm}
\includegraphics[width=0.20\textwidth]{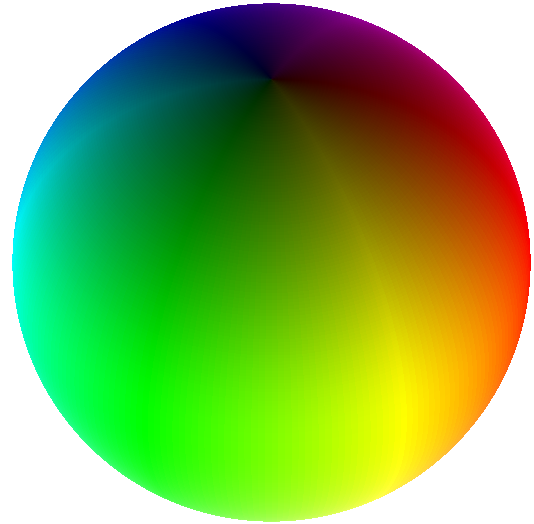}
Figure~6: Color schemes for phase plots and domain coloring on the
Riemann sphere
\end{center}
Somewhat surprisingly, the number of people using phase plots seems to
be quite small. The web site of Fran\c{c}ois Labelle \cite{Lab}
has a nice gallery of nontrivial pictures%
\footnote{Labelle justifies the sole use of phase by reasons
of clarity and aesthetics.},
including Euler's Gamma and Riemann's Zeta function.

Since the phase of a function occupies only one dimension of the
color space, there is plenty of room for depicting additional information.
It is recommended to encode this information by a gray scale, since
color (HUE) and brightness are {\em visually orthogonal}.
Figure~7 shows two such color schemes on the Riemann $w$-sphere.
\begin{center}\label{f.fig7}
\includegraphics[width=0.20\textwidth]{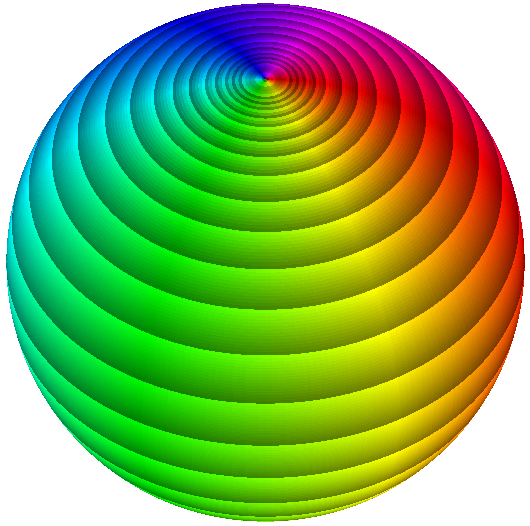}
\hspace{10mm}
\includegraphics[width=0.20\textwidth]{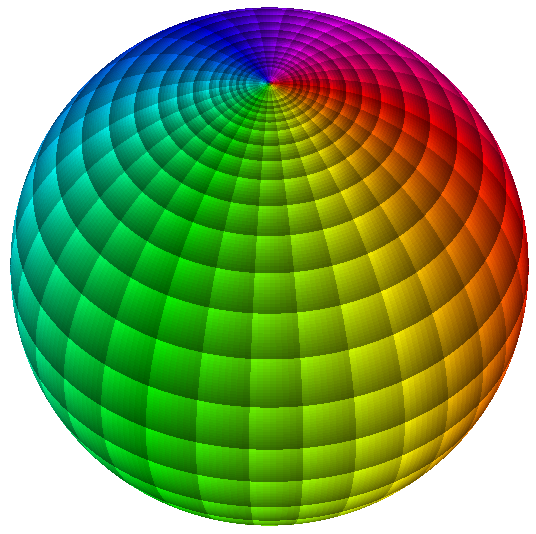}\\
Figure~7: Two color schemes involving sawtooth functions of gray
\end{center}
The left scheme is a combination of phase plots and standard domain coloring.
Here the
value $g$ of gray does not depend directly on $\log |w|$, but is a sawtooth
function thereof, $g(w) = \log |w|-\lfloor \log |w| \rfloor$.
This coloring works equally well, no matter in which range the values of
the function are located.

In the right scheme the gray value is the product of two sawtooth functions
depending on $\log |w|$ and $w/|w|$, respectively. The discontinuities of
this shading generate a logarithmically scaled polar grid.
{\em Pulling back} the coloring from the $w$-sphere to the $z$-domain of $f$ by
the mapping $w=f(z)$ resembles a {\em conformal grid mapping},
another well-known technique for depicting complex functions
(see Douglas Arnold \cite{Arn}). Note that pulling back a grid instead of
pushing it forward avoids multiple coverings.
Of course all coloring schemes can also be applied to functions on
Riemann surfaces.

For comparison, the figure below shows the four representations of
$f(z):=(z-1)/(z^2+z+1)$ in the square $|{\rm Re}\,z|\le 2$,
$|{\rm Im}\,z|\le 2$ corresponding to the color schemes of Figure~6 and
Figure~7, respectively.
\begin{center}\label{f.fig8}
\includegraphics[height=0.22\textwidth]{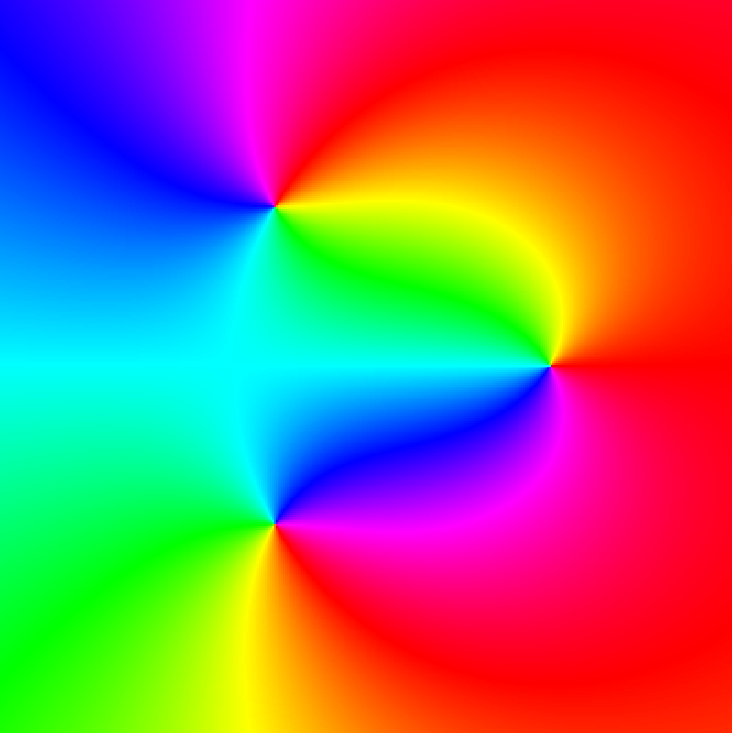}
\hfill
\includegraphics[height=0.22\textwidth]{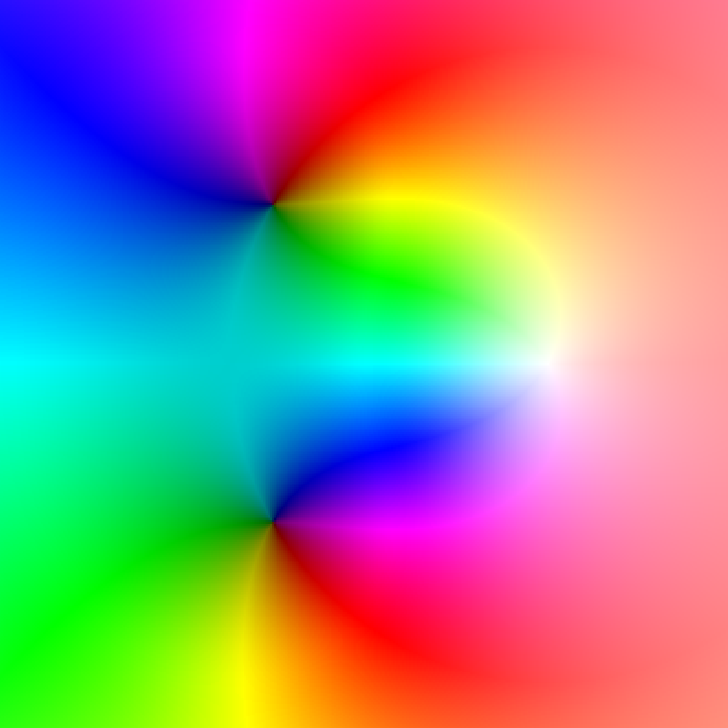}\\[8mm]
\includegraphics[height=0.22\textwidth]{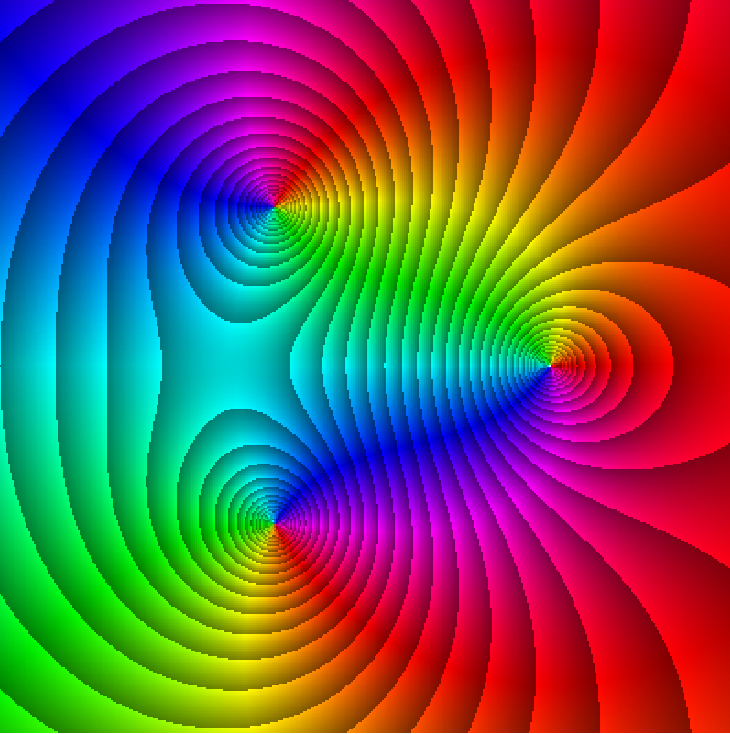}
\hfill
\includegraphics[height=0.22\textwidth]{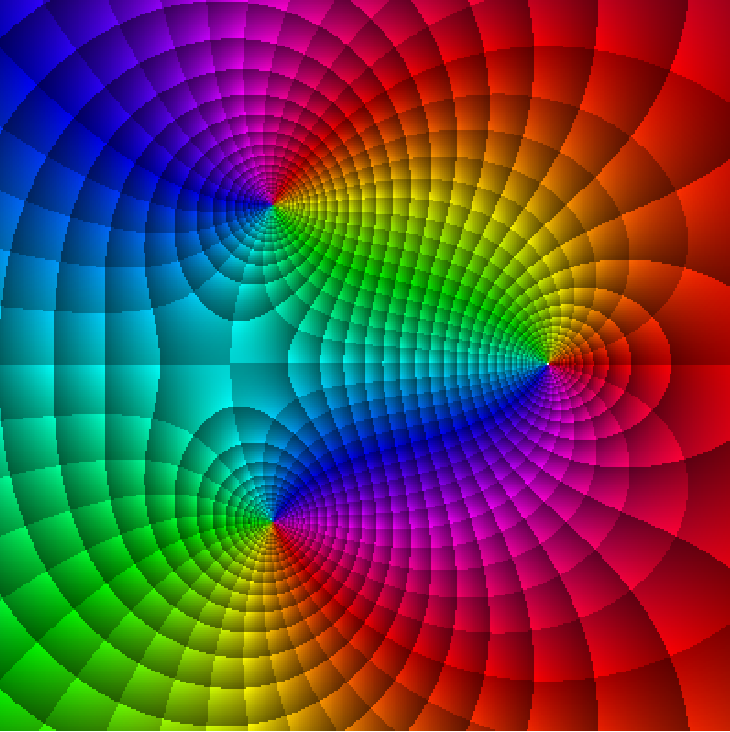}\\
Figure~8: Four representations of the function $f(z)=(z-1)/(z^2+z+1)$
\end{center}
Though these pictures (in particular the upper two) look quite similar,
which makes it simple to use them in parallel, the philosophy and the
mathematics behind them is quite different. We shall comment on this issue
in the final section.

\section*{The Phase Plot}

The phase of a complex function $f:D\rightarrow\widehat{\mathbb{C}}$
is defined on $D_0:=\{z\in D:\ f(z)\in \mathbb{C}^\times\}$, where
$\mathbb{C}^\times$ denotes the complex plane punctured at the
origin. Nevertheless we shall speak of phase plots
$P:D\rightarrow\mathbb{T},z\mapsto f(z)/|f(z)|$ on $D$, considering
those points where the phase is undefined as singularities. Recall
that $\mathbb{T}$ stands for the (colored) unit circle.

To begin with we remark that meromorphic functions are characterized
almost uniquely by their phase plot.

\begin{theorem}
If two non--zero meromorphic functions $f$ and $g$ on a connected
domain $D$ have the same phase, then $f$ is a positive scalar
multiple of $g$.
\end{theorem}

Proof.
Removing from $D$ all zeros and poles of $f$ and $g$ we get a connected
domain $D_0$. Since, by assumption, $f(z)/|f(z)|=g(z)/|g(z)|$ for all
$z\in D_0$, the function $f/g$ is holomorphic and real-valued in $D_0$, and
so it must be a (positive) constant.
\medskip

It is obvious that the result extends to the case where the phases of $f$ and
$g$ coincide merely on an open subset of $D$.

In order to check if two functions $f$ and $g$ with the same phase are equal,
it suffices to compare their values at a single point which is neither a zero
nor a pole. For purists there is also an intrinsic test which works with
phases alone: Assume that the non-constant meromorphic functions $f$ and $g$
have the same phase plot. Then it follows from the open mapping principle that
$f\not=g$ if and only if the phase plots of $f+c$ and $g+c$ are different
for one, and then for all, complex constants $c\not=0$.

\subsection*{Zeros and Poles}

Since the phases of zero and infinity are undefined, zeros and poles of a
function are singularities of its phase plot. What does the plot look like in a
neighborhood of such points?

If a meromorphic function $f$ has a zero of degree $n$ at $z_0$ it can be
represented as
\[
f(z)=(z-z_0)^n\,g(z),
\]
where $g$ is meromorphic and $g(z_0)\in \mathbb{C}^\times$. It follows that
the phase plot of $f$ close to $z_0$ looks like the phase plot of $z^n$ at
$0$, rotated by the angle $\mathrm{arg}\,g(z_0)$.
The same reasoning, with a negative integer $n$, applies to poles.
\begin{center}\label{f.fig9}
\includegraphics[width=0.4\textwidth]{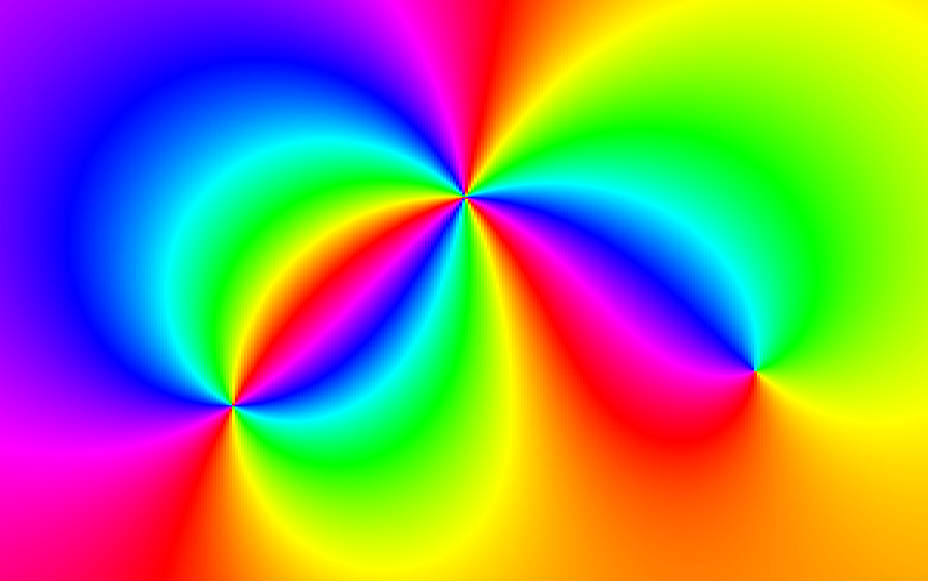} \\
Figure~9: A function with a simple zero, a double zero, and a triple pole
\end{center}
\medskip

Note that the colors are arranged in opposite orders for zeros and
poles. It is now clear that the phase plot does not only show
the location of zeros and poles but also reveals their {\em multiplicity}.

A useful tool for locating zeros is the {\em argument principle}. In order to
formulate it in the context of phase plots we translate the definition of
winding number into the language of colors:
Let $\gamma:\mathbb{T}\rightarrow D_0$ be a closed oriented path in the domain
$D_0$ of a phase plot $P:D_0\rightarrow\mathbb{T}$. Then the usual
winding number of the mapping $P\circ\gamma:\mathbb{T}\rightarrow
\mathbb{T}$ is called the {\em chromatic number} of $\gamma$ with respect to the
phase plot $P$ and is denoted by $\mathrm{chrom}_P\gamma$ or simply by
$\mathrm{chrom}\,\gamma$.

Less formally, the chromatic number counts how many times the color of
the point $\gamma(t)$ moves around the complete color circle when
$\gamma(t)$ traverses $\gamma$ once in positive direction.

Now the argument principle can be rephrased as follows:
{\em
Let $D$ be a Jordan domain with positively oriented boundary $\partial D$
and assume that $f$ is meromorphic in a neighborhood of $D$.
If $f$ has $n$ zeros and $p$ poles in $D$ (counted with multiplicity), and none
of them lies on $\partial D$, then}
\vspace{-2mm}
\[
n-p = \mathrm{chrom}\,\partial D.
\]
\begin{center}\label{f.fig10}
\includegraphics[width=0.4\textwidth]{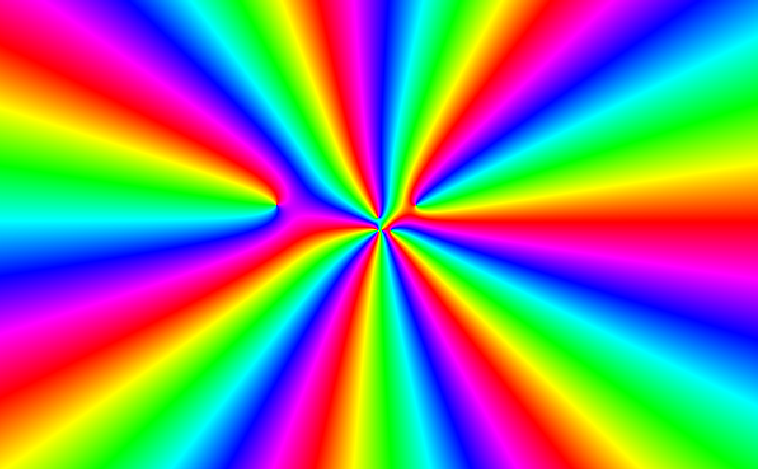} \\
Figure~10: This function has no poles. How many zeros are in the
displayed rectangle?
\end{center}

Looking at Figure~10 in search of zeros immediately brings forth
new questions, for example:
Where do the isochromatic lines end up? Can they connect two
zeros? If so, do these lines have a special meaning? What about ``basins
of attraction''? Is there always a natural (cyclic) ordering of zeros?
What can be said about the global structure of phase plots? We shall
return to these issues later.

\subsection*{The Logarithmic Derivative}

Along the {\em isochromatic lines} of a phase plot the argument of $f$ is
constant. The Cauchy-Riemann equations for any continuous branch of the
logarithm $\log f = \ln |f| + \mathrm{i}\,\arg f$
imply that these lines are orthogonal to the level lines of $|f|$,
i.e.\ the isochromatic lines are parallel to the gradient of $|f|$.
According to the chosen color scheme, we have red on the right and green on the
left when walking on a yellow line in ascending direction.

To go a little beyond this qualitative result, we denote by $s$ the unit
vector parallel to the gradient of $|f|$ and set $n:=\mathrm{i} s$.
With $\varphi := \arg f$ and $\psi:=\log |f|$ the Cauchy-Riemann equations
for $\log f$ imply that the directional derivatives of $\varphi$ and $\psi$
satisfy
\[
\partial_s \psi = \partial_n \varphi > 0, \quad
\partial_n \psi = -\partial_s \varphi = 0,
\]
at all points $z$ of the phase plot where $f(z)\not=0$ and $f'(z)\not=0$.
Since the absolute value of $\partial_n \varphi$
measures the density of the isochromatic lines, we can visually estimate
the growth of $\log |f|$ along these lines from their density.
Because the phase plot delivers no information on the absolute value,
this does not say much about the growth of $|f|$.
But taking into account the second Cauchy-Riemann equation and
\[
|(\log f)'|^2 = (\partial_n \varphi)^2  + (\partial_s \varphi)^2,
\]
we obtain the correct interpretation of the density $\partial_n \varphi$:
it is the modulus of the {\em logarithmic derivative},
\begin{equation} \label{e.dens}
\partial_n \varphi = \left|{f'}/{f}\right|.
\end{equation}
So, finally, we need not worry about branches of the logarithm.
It is worth mentioning that $\partial_n \varphi (z)$ behaves asymptotically
like $k/|z-z_0|$ if $z$ approaches a zero or pole of order $k$ at $z_0$.

But this is not yet the end of the story. What about zeros of $f'$\,?
Equation (\ref{e.dens}) indicates that something should be visible
in the phase plot. Indeed, points $z_0$ where $f'(z_0)=0$ and $f(z_0)\not=0$
are ``{\em color saddles}'', i.e.\, intersections of isochromatic lines.

If $f'$ has a zero of order $k$ at $z_0$, then $z\mapsto f(z)-f(z_0)$ has a
zero of order $k+1$ at $z_0$. Consequently $f$ can be represented as
\[
f(z) = f(z_0) + (z-z_0)^{k+1}\,g(z)
\]
where $g(z_0)\not=0$. It follows that $f(z)$ travels $k+1$ times around
$f(z_0)$ when $z$ moves once around $z_0$ along a small circle.
In conjunction with $f(z_0)\not=0$ this can be used to show
that there are exactly $2k+2$ isochromatic lines emanating from $z_0$
where the phase of $f$ is equal to the phase of $f(z_0)$.
Alternatively, one can also think of $k+1$ smooth isochromatic lines
intersecting each other at $z_0$.
\vspace*{-4mm}
\begin{center}\label{f.fig13}
\includegraphics[width=0.22\textwidth]{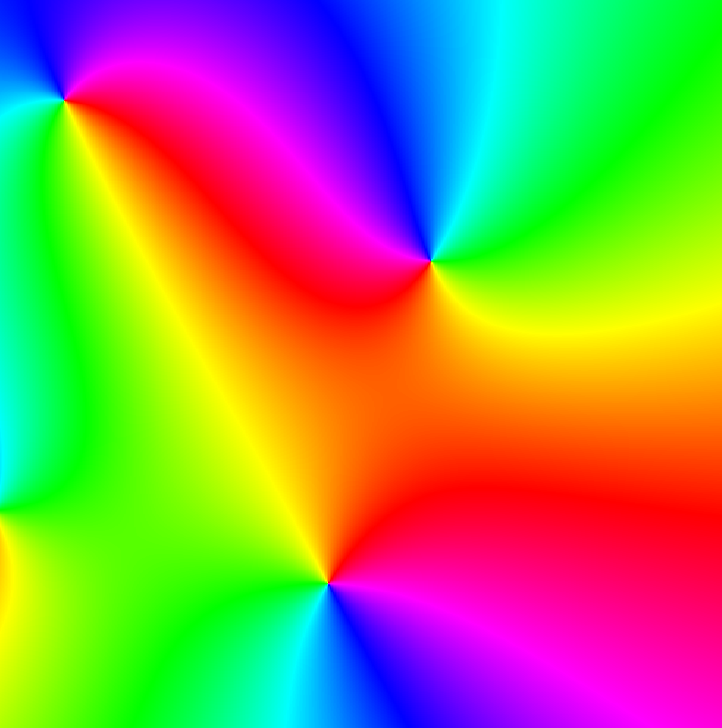} \hfill
\includegraphics[width=0.22\textwidth]{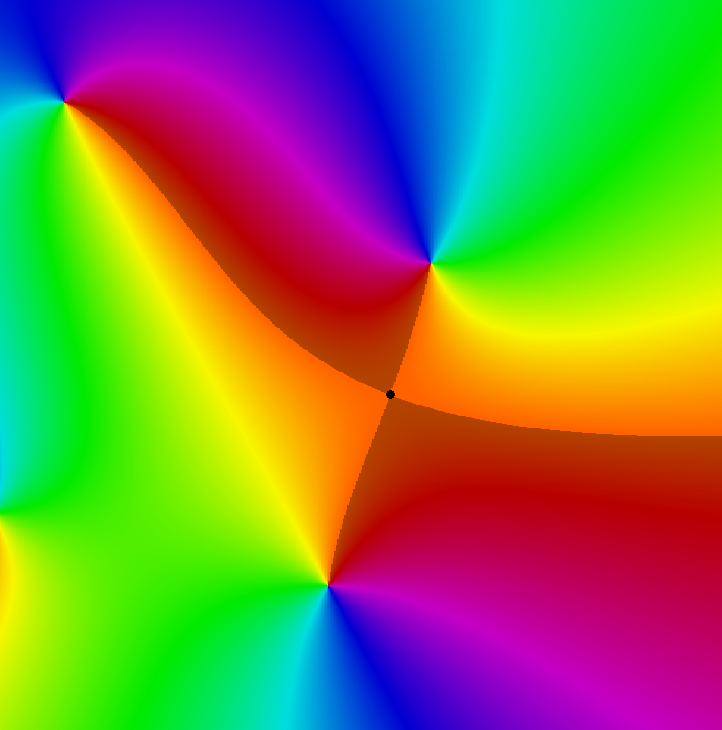} \\
Figure~11: Zeros of $f'$ are color saddles
\end{center}
Color saddles appear as diffuse spots like in the left picture of Figure~11.
To locate them
precisely it is helpful to modify the color scheme by superimposing a gray
component
which has a jump at some point $t$ of the unit circle. If $t:=f(z_0)/|f(z_0)|$
is chosen, then the phase plot shows a sharp saddle at the zero
$z_0$ of $f'$ as in the right picture.

\subsection*{Essential Singularities}

Have you ever seen an essential singularity?
Here is the picture which usually illustrates this situation.

\begin{center}\label{f.fig11}
\includegraphics[width=0.35\textwidth]{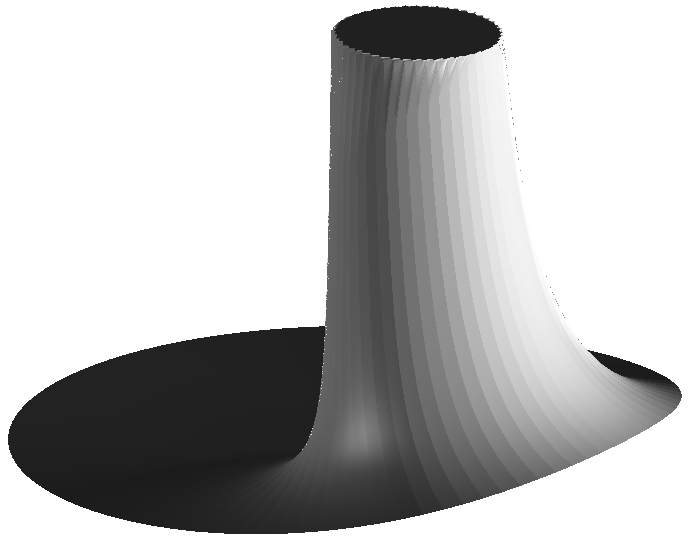} \\
Figure~12: The analytic landscape of $f(z)=\mathrm{e}^{1/z}$
\end{center}
Despite the massive tower this is not very impressive, and with
regard to the Casorati-Weierstrass Theorem or the Great Picard
Theorem one would expect something much wilder.
Why does the analytic landscape not reflect this behavior?
For the example the answer is easy: the function has a tame modulus,
every contour line is a single circle through the origin.
Now look at the phase plot in Figure~13:
\begin{center}\label{f.fig12}
\includegraphics[width=0.30\textwidth]{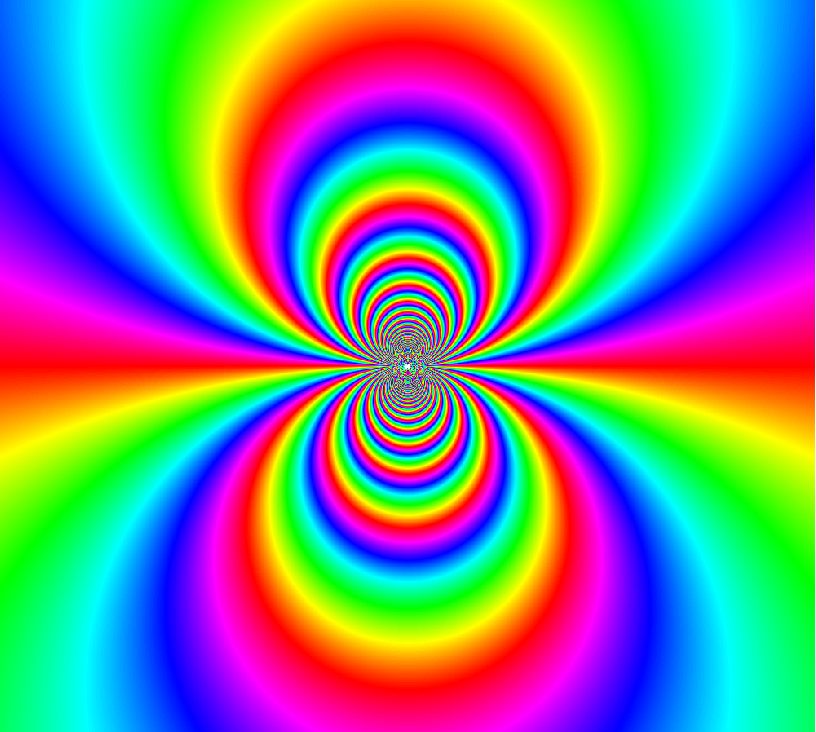} \\
Figure~13: A phase plot depicting the essential singularity of
$f(z)=\mathrm{e}^{1/z}$
\end{center}
But must there not be a symmetry between modulus and phase? In
fact not. There is such a symmetry of modulus and {\em argument}
(for non-vanishing functions), but phase plots depict the phase and
not the argument -- and this makes a difference.

So much for the example, but what about the general case? Perhaps there
are also functions which conceal their essential singularities in the
phase plot?

In order to show that this cannot happen, we assume that
$f: D\rightarrow\mathbb{C}$ is analytic and has an essential singularity
at $z_0$.

By the Great Picard Theorem, there exists a color $c\in\mathbb{T}$ such
that any punctured neighborhood $U$ of $z_0$ contains infinitely many points
$z_k\in U$ with $f(z_k)=c$. Moreover, the set of zeros of $f'$ in $D$
is at most countable, and hence we can choose $c$ such that
$c \not= f(z)/|f(z)|$ for all zeros $z$ of $f'$.

Since $|c|=1$, all isochromatic lines of the phase plot through the points
$z_k$ have the color $c$. As was shown in the preceding section, the
modulus of $f$ is strictly monotone along these lines. Now it is not hard
to see that two distinct points $z_k$ cannot lie on the same isochromatic
line, because these lines can meet each other only at a zero of $f'$,
which has been excluded by the special choice of $c$.

Consequently any neighborhood of an essential singularity contains a countable
set of pairwise disjoint isochromatic lines with color $c$.
Combining this observation with the characterization of phase plots near
poles and removable singularities we obtain the following result.

\begin{theorem}\label{t.ess}
An isolated singularity $z_0$ of an analytic function $f$ is an essential
singularity if and only if any neighborhood of $z_0$ intersects infinitely
many isochromatic lines of the phase plot with one and the same color.
\end{theorem}

Note that a related result does not hold for the argument, since then,
in general, the values of $\arg f(z_k)$ are different. For example,
any two isochromatic lines of the function $f(z)=\exp(1/z)$ have a
different argument.

\subsection*{Periodic Functions}

Obviously, the phase of a periodic function is periodic, but what about the
converse?

Though there are only {\em two} classes (simply and doubly periodic) of
nonconstant periodic meromorphic functions on $\mathbb{C}$, we can observe
{\em three} different types of periodic phase plots.

\begin{center}\label{f.fig14}
\includegraphics[width=0.40\textwidth]{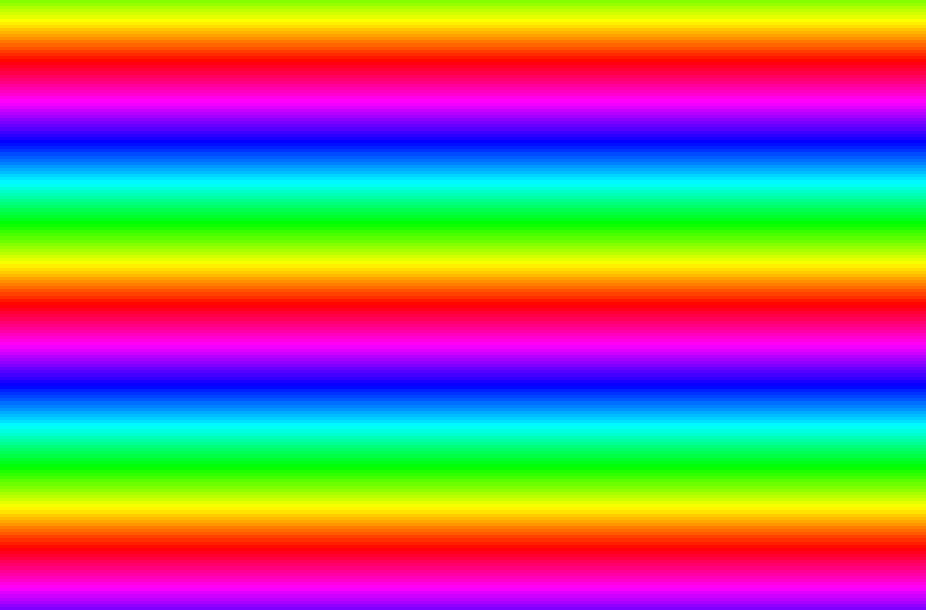} \\
Figure~14: Phase plot of $f(z)=\mathrm{e}^z$
\end{center}

\begin{center}\label{f.fig15}
\includegraphics[width=0.40\textwidth]{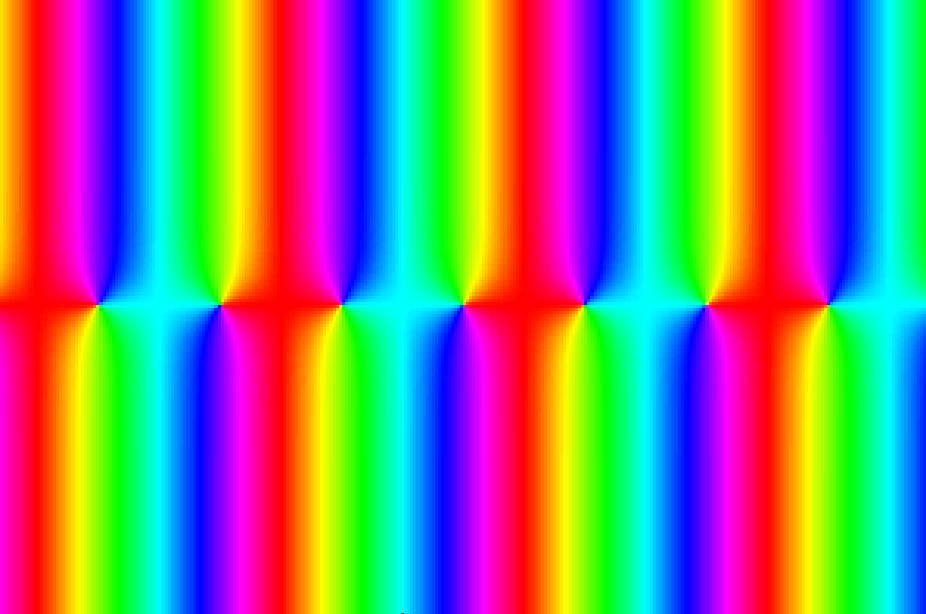} \\
Figure~15: Phase plot of $f(z)=\sin z$
\end{center}

\begin{center}\label{f.fig16}
\includegraphics[width=0.40\textwidth]{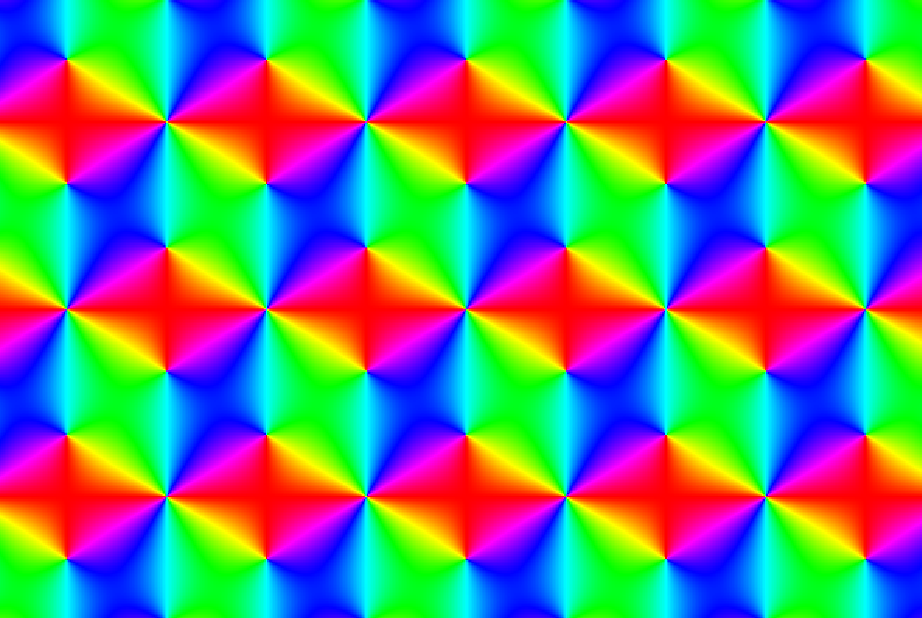} \\
Figure~16: Phase plot of a Weierstrass $\wp$-function
\end{center}

``Striped'' phase plots like in Figure~14 always depict exponential
functions $f(z) = \mathrm{e}^{az+b}$ with $a\not=0$. Functions with simply
$p$--periodic phase need not be periodic, but have the more general form
$\mathrm{e}^{\alpha z/p}\,g(z)$ with $\alpha\in\mathbb{R}$ and a
$p$--periodic function $g$.
Somewhat surprisingly, doubly periodic phase plots indeed {\em always}
represent elliptic functions.

The first result basically follows from the fact that the function
$\arg f$ is harmonic and has parallel straight contour lines, which implies
that $\arg f(x+{\mathrm{i}} y) = \alpha x + \beta y + \gamma$. Since
$\log |f|$ is conjugate harmonic to $\arg f$, it necessarily has the form
$\log |f(x+{\mathrm{i}} y)| = - \alpha y + \beta x + \delta$.

If the phase of $f$ is $p$-periodic, then we have
\begin{equation*}\label{e.hhh}
h(z):=\frac{f(z+p)}{f(z)} = \frac{|f(z+p)|}{|f(z)|} \in \mathbb{R}_+,
\end{equation*}
and since $h$ is meromorphic on $\mathbb{C}$, it must be a positive constant
$\mathrm{e}^\alpha$. Now it follows easily that
$g(z):=f(z)\cdot\mathrm{e}^{-\alpha z/p}$ is periodic with period $p$.

Finally, if $p_1$ and $p_2$ are periods of $f/|f|$ with
$p_1/p_2\notin\mathbb{R}$, then
there exist $\alpha_1,\alpha_2\in\mathbb{R}$ such that
$f(z+p_j)=\mathrm{e}^{\alpha_j}\,f(z)$.
The meromorphic function $g$ defined by $g(z):=f'(z)/f(z)$ has only simple poles
and zeros. Integration of $g = (\log f)'$ along a (straight) line
from $z_0$ to $z_0+p_j$ which contains no pole of $g$ yields that
\[
\alpha_j = \int_{z_0}^{z_0+p_j} g(z)\,dz.
\]
Evaluating now the area integral $\iint_\Omega g\,dx\,dy$ over the
parallelogram $\Omega$ with vertices at $0,p_1,p_2,p_1+p_2$ by two
different iterated integrals,
we obtain $\alpha_2\,p_1 = \alpha_1\,p_2$.
Since $\alpha_1,\alpha_2\in\mathbb{R}$ and $p_1/p_2\notin\mathbb{R}$ this
implies that $\alpha_j=0$.

\subsection*{Partial Sums of Power Series}

Figure~17 shows a strange image which, in similar form, occurred in an
experiment. Since it looks so special, one could attribute it
to a programming error.
A moment's thought reveals what is going on here, at least at an intuitive
level. This example demonstrates again that looking at phase plots can
immediately provoke new questions.

\begin{center}\label{f.fig17}
\includegraphics[width=0.35\textwidth]{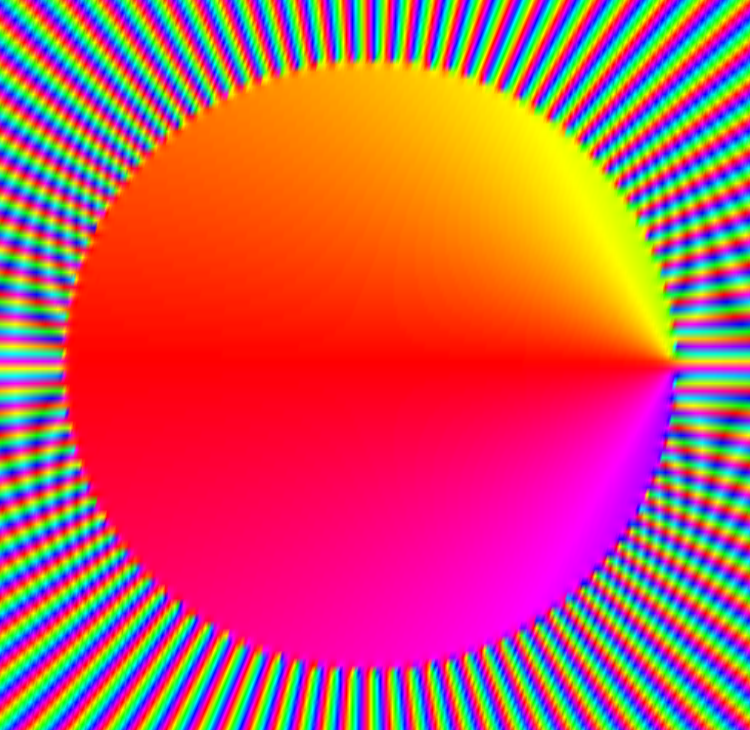} \\
Figure~17: A Taylor polynomial of $f(z)=1/(1-z)$
\end{center}

Indeed the figure illustrates a rigorous result (see Titchmarsh \cite{Tit},
Section 7.8) which was proven by Robert Jentzsch in 1914:
\medskip

{\em If a power series $a_0+a_1z+a_2z^2+\ldots$ has a positive finite
convergence radius $R$, then the zeros of its partial sums cluster at every
point $z$ with $|z|=R$.}
\medskip

The reader interested in the life and personality of Robert Jentzsch is
referred to the recent paper \cite{DuHeKha} by Peter Duren, Anne-Katrin Herbig,
and Dmitry Khavinson.

\section*{Boundary Value Problems}

Experimenting with phase plots raises a number of new questions.
One such problem is to find a criterion for deciding which color images
are {\em analytic phase plots}, i.e., phase plots of analytic functions.

Since phase plots are painted with the restricted palette of
{\em saturated} colors from the color circle, Leonardo's Mona Lisa
will certainly never appear.
But for analytic phase plots there are much stronger restrictions:
By the uniqueness theorem for harmonic functions an arbitrarily small
open piece determines the plot entirely.

So let us pose the question a little differently: What are appropriate
data which can be prescribed to construct an analytic phase plot, say,
in a Jordan domain $D$\,? Can we start, for instance, with  given colors
on the boundary $\partial D$?
If so, can the boundary colors be prescribed arbitrarily or are they
subject to constraints?

In order to state these questions more precisely we introduce the concept of
a {\em colored set} $K_C$, which is a subset $K$ of the complex plane
together with a mapping $C: K\rightarrow\mathbb{T}$. Any such mapping is
referred to as a {\em coloring} of $K$.

For simplicity we consider here only the following setting of boundary value
problems for phase plots with continuous colorings:
\medskip\smallskip

{\em
Let $D$ be a Jordan domain and let $B$ be a continuous coloring
of its boundary $\partial D$.
Find all continuous colorings $C$ of $\overline{D}$ such that the
restriction of $C$ to $\partial D$ coincides with $B$ and the restriction
of $C$ to $D$ is the phase plot of an analytic function $f$ in $D$.
}
\medskip\smallskip

If such a coloring $C$ exists, we say that the coloring $B$ {\em admits a
continuous analytic extension to} $\overline{D}$.

The restriction to continuous colorings automatically excludes zeros
of $f$ in $D$. It does, however, {\em not} imply that $f$ must extend
continuously onto $\overline{D}$ -- and in fact it is essential
not to require the continuity of $f$ on $\overline{D}$ in order to get a
nice result.

\begin{theorem}
\label{t.rhp}
Let $D$ be a Jordan domain with a continuous coloring
$B$ of its boundary $\partial D$. Then $B$ admits a continuous
analytic extension to $\overline{D}$ if and only if the chromatic
number of $B$ is zero. If such an extension exists, then it is unique.
\end{theorem}

Proof.
If $C:\overline{D}\rightarrow\mathbb{T}$ is a continuous coloring, then
a simple homotopy argument (contract $\partial D$ inside $D$ to a point)
shows that the chromatic number of its restriction to $\partial D$ must vanish.

Conversely,
any continuous coloring $B$ of $\partial D$ with chromatic number zero
can be represented as $B=\mathrm{e}^{\mathrm{i}\varphi}$ with a continuous
function
$\varphi: \partial D\rightarrow\mathbb{R}$. This function admits a
unique continuous harmonic extension $\Phi$ to $\overline{D}$.
If $\Psi$ denotes a harmonic conjugate of $\Phi$, then
$f=\mathrm{e}^{\mathrm{i}\Phi-\Psi}$ is
analytic in $D$. Its phase $C:=\mathrm{e}^{\mathrm{i}\Phi}$ is continuous on
$\overline{D}$
and coincides with $B$ on $\partial D$.
\medskip

Theorem~\ref{t.rhp} {\em parametrizes} analytic phase plots which
extend continuously on $\overline D$ by their boundary colorings.
This result can be generalized to phase plots which are continuous on
$\overline{D}$ with the exception of finitely many singularities of zero
or pole type in $D$. Admitting now boundary colorings $B$ with arbitrary
color index we get the following result:

{\em For any finite collection of given zeros with orders $n_1,\ldots,n_j$
and poles of orders $p_1,\ldots,p_k$ the boundary value problem for
meromorphic phase plots with prescribed singularities has a unique
solution if and only if the (continuous) boundary coloring $B$ satisfies}
\[
\mathrm{chrom}\,B = n_1+\ldots+n_j-p_1-\ldots -p_k.
\]


\section*{The Riemann Zeta Function}

After these preparations we are ready to pay a visit to ``Zeta'', the
mother of all analytic functions. Here is a phase plot in the square
$-40\le{\rm Re}\,z\le 10,$ $-2\le{\rm Im}\,z\le 48$.
\begin{center}\label{f.fig21}
\includegraphics[width=0.45\textwidth]{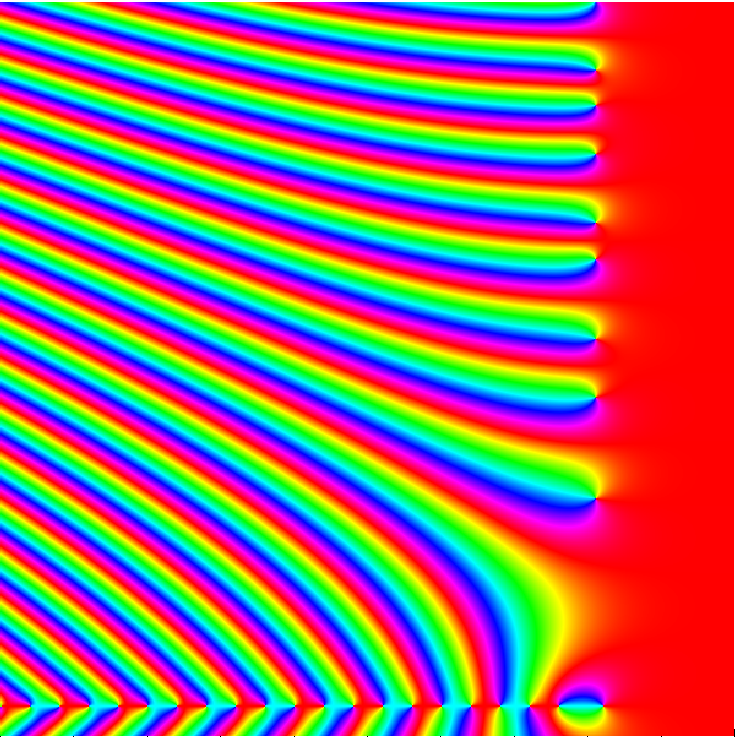}\\
Figure~18: The Riemann Zeta function
\end{center}
We see the pole at $z=1$, the {\em trivial zeros} at the points
$-2,-4,-6,\ldots $ and several zeros on the {\em critical line}
$\mathrm{Re}\,z=1/2$. Also we observe that the isochromatic lines are
quite regularly distributed in the left half plane.

Saying that Zeta is the mother of all functions alludes to its
{\em universality}.
Our starting point is the following strong version of Voronin's Universality
Theorem due to Bagchi \cite{Bag} (see also Karatsuba and Voronin \cite{KarVor},
Steuding \cite{Steu}):
\medskip\smallskip

{\em
Let $D$ be a Jordan domain such that $\overline{D}$ is contained in the strip
\[
R:=\{z\in \mathbb{C}: 1/2 < \mathrm{Re}\, z < 1 \},
\]
and let $f$ be any function which is analytic in $D$, continuous on
$\overline{D}$, and has no zeros in $\overline{D}$. Then $f$ can be uniformly
approximated on $\overline{D}$ by vertical shifts of Zeta,
$\zeta_t(z):=\zeta(z+\mathrm{i} t)$ with $t\in\mathbb{R}$.
}
\medskip\smallskip

Recall that a continuously colored Jordan curve $J_C$ is a continuous mapping
$C: J\rightarrow\mathbb{T}$ from a simple closed curve $J$ into the color circle
$\mathbb{T}$.\\[-3mm]

\begin{minipage}{0.2\textwidth}
\begin{center}\label{f.fig20a}
\includegraphics[width=0.7\textwidth]{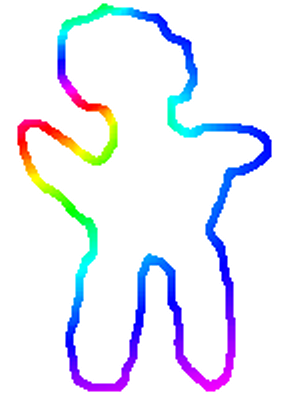}\\
Figure~19: A representative of a string
\end{center}
\end{minipage}
\hfill
\begin{minipage}{0.25\textwidth}
A {\em string} $S$ is an equivalence class of all such colored curves
with respect to rigid motions of the plane.
Like colored Jordan curves, strings fall into different classes according to
their {\em chromatic number}. Figure~19 depicts a representative of a
string with chromatic number one.
\end{minipage}
\vspace*{3mm}

We say that a string $S$ {\em lives in a domain} $D$ if it can be
represented by a colored Jordan curve $J_C$ with $J\subset D$. A string
can {\em hide itself} in a phase plot $P: D\rightarrow\mathbb{T}$,
if, for every $\varepsilon>0$, it has a representative $J_C$ such that
$J\subset D$ and
\[
\max_{z\in J} |C(z)-P(z)|<\varepsilon.
\]
In less technical terms, a string can hide itself if it can move to
a place where it is invisible since it blends in almost perfectly
with the background.

In conjunction with Theorem~\ref{t.rhp} the following universality result
for the phase plot of the Riemann Zeta function can easily be derived
from Voronin's theorem.

\begin{theorem}
\label{t.uni}
Let $S$ be a string which lives in the strip $R$. Then $S$ can hide
itself in the phase plot of the Riemann Zeta function on $R$ if it
has chromatic number zero.
\end{theorem}

In view of the extreme richness of Jordan curves and colorings this result
is a real miracle.
The three pictures below show phase plots of Zeta in the critical strip.
The regions with saturated colors belong to $R$.
The rightmost figure depicts the domain considered on p.\,342 of Conrey's
paper \cite{Con}.

\begin{center}\label{f.fig19}
\includegraphics[height=0.4\textwidth]{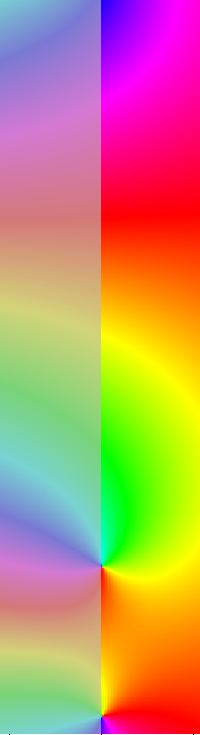}\hfill
\includegraphics[height=0.4\textwidth]{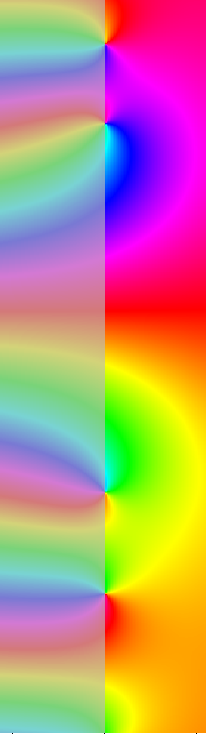}\hfill
\includegraphics[height=0.4\textwidth]{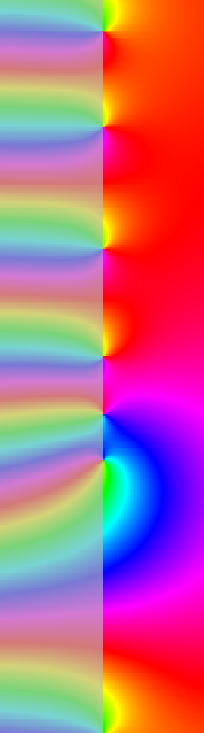}\\
Figure~20: The Riemann Zeta function at $\mathrm{Im}\,z=171$,
$8230$ and $121415$
\end{center}

What about the converse of Theorem~\ref{t.uni}\,? If there existed
strings with {\em nonzero chromatic number} which can hide themselves
in the strip $R$, their potential hiding-places must be Jordan curves
with non-vanishing chromatic number in the phase plot.
By the argument principle, this would imply that Zeta has zeros in $R$.
If we assume this, for a moment, then such strings indeed exist:
They are perfectly hidden and wind themselves once around such a zero.
So the converse of Theorem~\ref{t.uni} holds if and only if $R$ contains
no zeros of Zeta, which is known to be equivalent to the {\em Riemann
hypothesis} (see Conrey \cite{Con}, Edwards \cite{Edw}).

\section*{Phase Flow and Diagrams}

Mathematical creativity is based on the interplay of problem
posing and problem solving, and it is my belief that the former
is even more important than the latter: often the key to solving a
problem lies in asking the right questions.

Illustrations have a high density of information and stimulate imagination.
Looking
at pictures helps in getting an intuitive understanding of mathematical objects
and finding interesting questions, which then can be investigated using rigorous
mathematical techniques.

This section intends to demonstrate how phase plots can produce novel ideas.
The material presented here is the protocol of a self-experiment which has been
carried out by the author in order to check the creative potential of phase
plots.

Let us start by looking at Figure~1 again. It depicts the phase of
a {\em finite Blaschke product}, which is a function of the form
\[
f(z) = c\,\prod_{k=1}^n \frac{z-z_k}{1-\overline{z_k}z},\quad z\in\mathbb{D},
\]
with $|z_k|<1$ and $|c|=1$. Blaschke products are fundamental building blocks
of analytic functions in the unit disc and have the property
$|f(z)|=1$ for all $z\in\mathbb{T}$. The function shown in Figure~1 has
$81$ zeros $z_k$ in the unit disk.

Looking at Figure~1 for a while leaves the impression of a cyclic ordering
of the zeros. Let us test this with another example having only five zeros
(Figure~21, left).
\begin{center}
\includegraphics[width=0.22\textwidth]{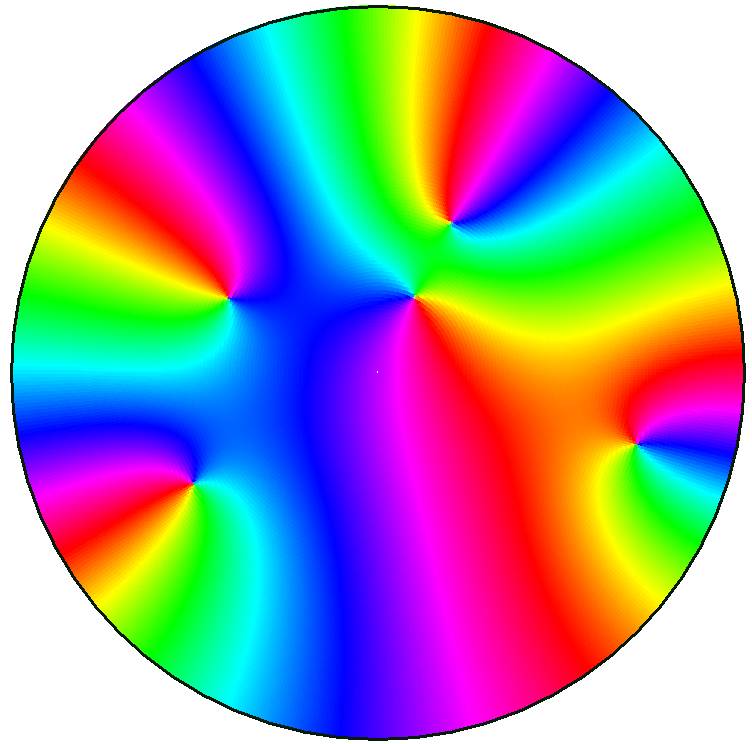}
\includegraphics[width=0.22\textwidth]{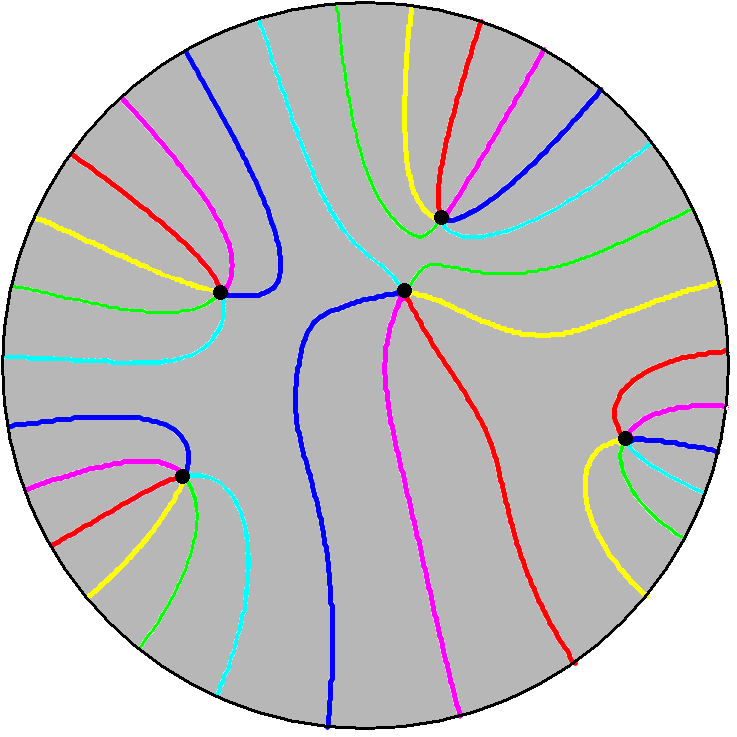} \\
Figure~21: The phase plot of a Blaschke product with five zeros
\end{center}
The left picture seems to confirm the expectation: if we focus attention
to the yellow color, any of these lines connects a zero with a
certain point on the boundary, thus inducing a cyclic ordering.

However, looking only at one specific color is misleading. Choosing
another, for instance blue, can result in a different ordering.
So what is going on here? More precisely:
What is the {\em global structure} of the phase plot of a Blaschke product?
This could be a good question.

An appropriate mathematical framework to develop this idea is the theory of
{\em dynamical systems}. We here only sketch the basic facts; for details
see \cite{Wgt2}.

With any meromorphic function $f$ in a domain $D$ we associate the
dynamical system
\begin{equation} \label{e.dynsys}
\dot{z}  = g(z) := \frac{f(z)\,\overline{f'(z)}}{|f(z)|^2+|f'(z)|^2}.
\end{equation}
The function $g$ on the right-hand side of (\ref{e.dynsys})
extends from $D_0$ to a smooth function on $D$.
This system induces a {\em flow} $\Phi$ on $D$,
which we designate as the {\em phase flow} of $f$.

The fixed points of (\ref{e.dynsys}) are the zeros of $f$ (repelling), the
poles of $f$ (attracting), and the zeros of $f'$ (saddles). The remaining
orbits are the components of the isochromatic lines of the phase plot of $f$
when the fixed points are removed. Thus the orbits of the phase flow
endow the phase plot with an additional structure and convert it into a
{\em phase diagram}.

Intuitively, the phase flow $\Phi$ transports a colored substance (``phase'')
from the zeros to the poles and to the boundary of the domain.
\begin{center}
\includegraphics[width=0.22\textwidth]{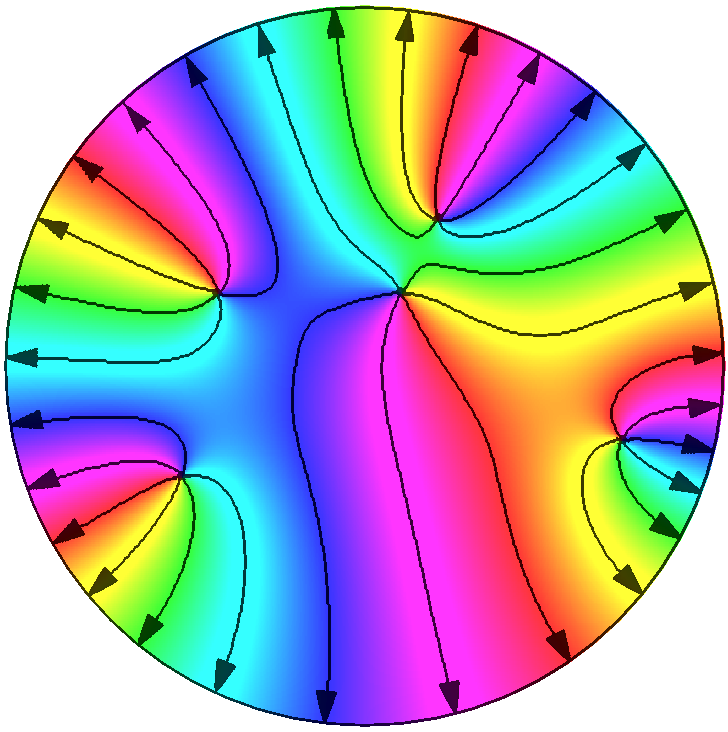}\hfill
\includegraphics[width=0.22\textwidth]{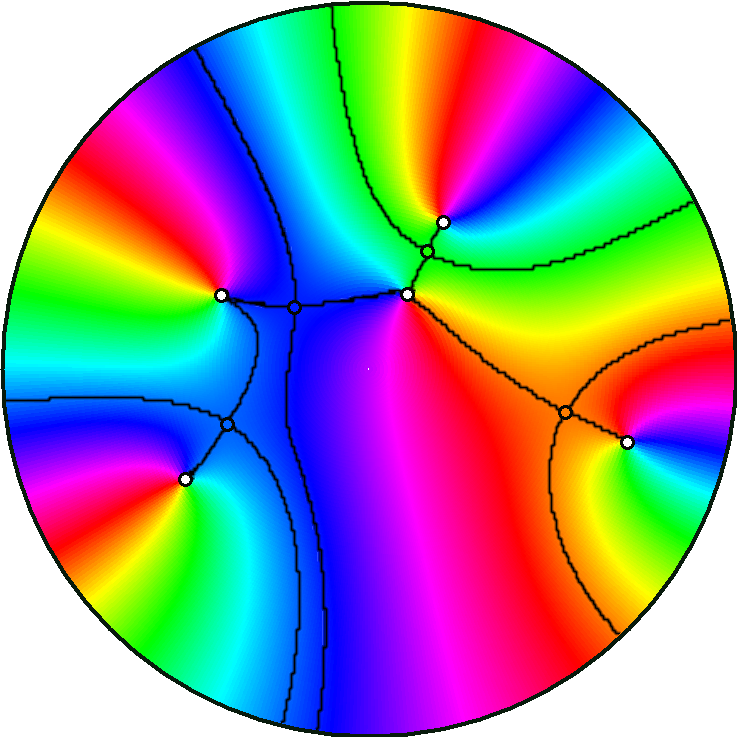}\\
Figure~22: Phase transport to the boundary, zeros of $f$ and $f'$ with
invariant manifolds
\end{center}
The left part of Figure~23 illustrates how ``phase'' of measure $2\pi$ emerging
from the zero in the highlighted domain is transported along the orbits of
$\Phi$ until it is finally deposited along (parts of) the boundary.

In the general case, where $f:D\rightarrow\widehat{\mathbb{C}}$ is meromorphic
on a domain $D$ and $G\subset D$ is a Jordan domain with boundary $J$ in $D_0$,
the phase flow of any zero (pole) of $f$ in $G$ generates a (signed) {\em
measure} on $J$.
The result is a {\em quantitative version of the argument principle} which
tells us in which way the phase of the zeros (poles) is distributed along
$J$ (see Figure~23, right).
\begin{center}
\includegraphics[height=0.22\textwidth]{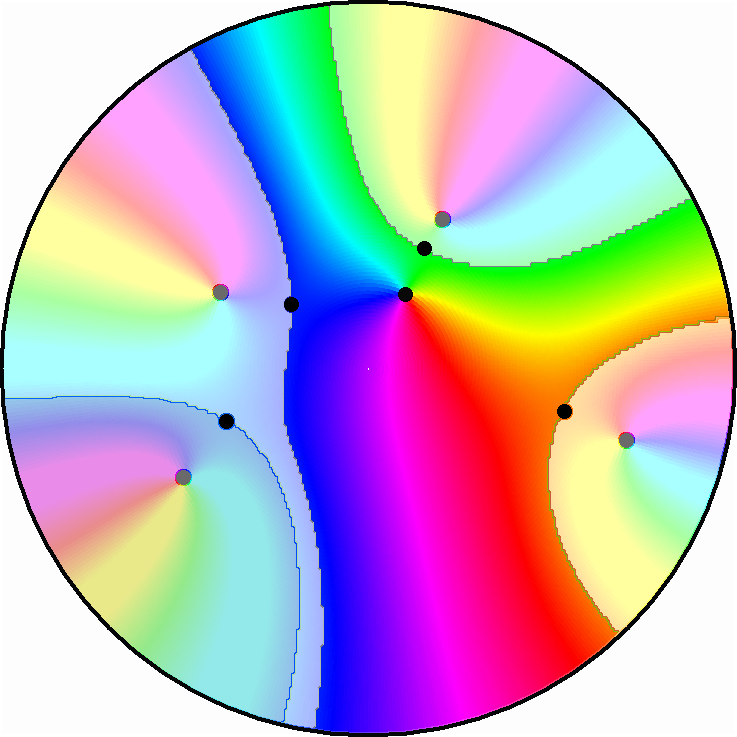} \hfill
\includegraphics[height=0.22\textwidth]{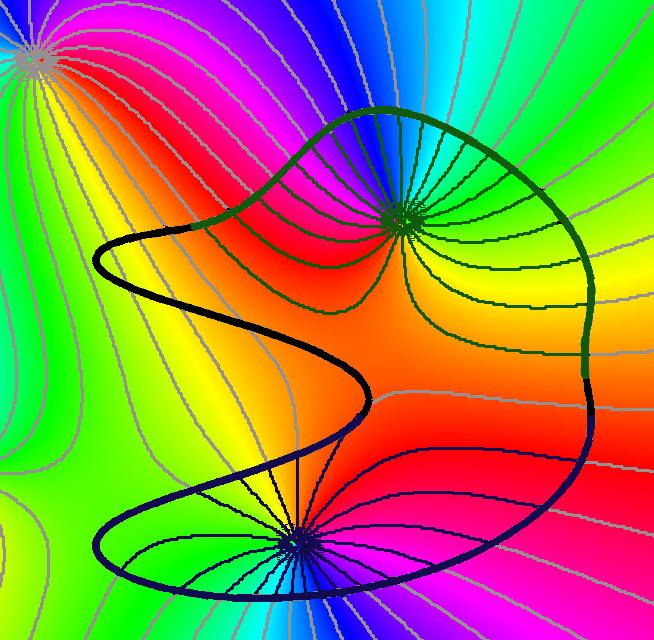} \\
Figure~23: The phase flow and the argument principle
\end{center}

The question about the structure of phase plots of Blaschke products can now 
be rephrased in the setting of dynamical systems:
What are the {\em basins of attraction} of the zeros of $f$ with respect
to the (reversed) phase flow?

The key for solving this problem is given by the {\em invariant manifolds}
of the saddle points, i.e., the points $a_j\in \mathbb{D}$ where $f'(a_j)=0$
and $f(a_j)\not=0$.
\begin{center}
\includegraphics[width=0.22\textwidth]{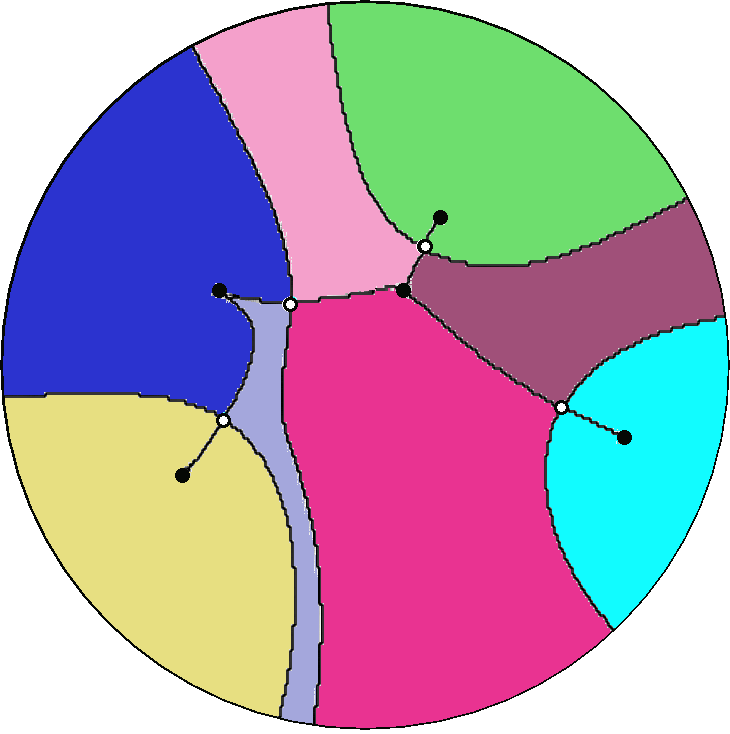}\hfill
\includegraphics[width=0.22\textwidth]{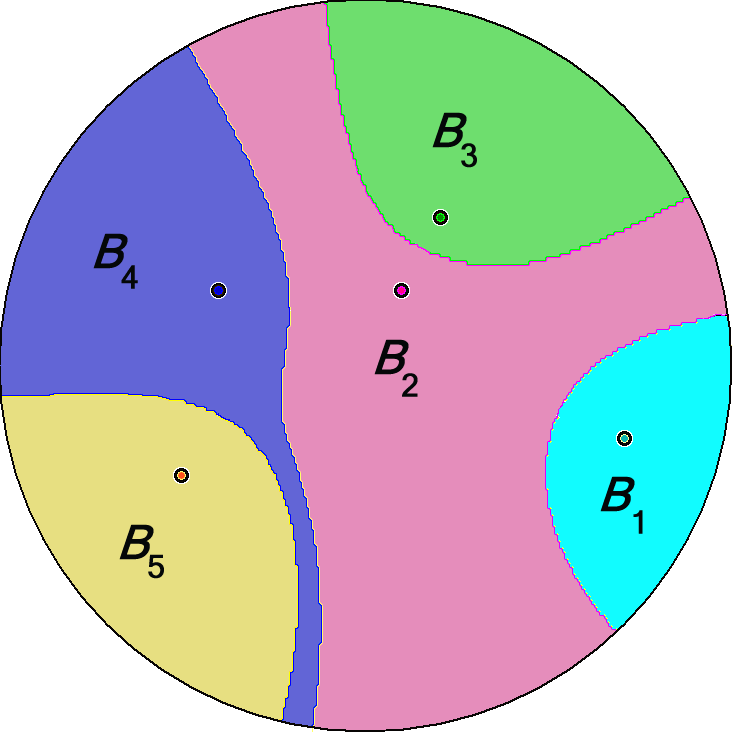} \\
Figure~24: Invariant manifolds of the saddles and basins of attraction of
the zeros
\end{center}
Removing all {\em unstable manifolds} of the points $a_j$ from $\mathbb{D}$
results in an open set $B$, which is the union of connected components $B_j$.
Any component $B_j$ contains exactly one zero $b_j$ of $f$, where multiple
zeros are counted only once.

The intersection of every set $\overline{B_j}$ with $\mathbb{T}$ is not
empty and consists of a finite number of arcs $A_{ji}$.
The complete set of these arcs covers the unit circle and two arcs
are either disjoint or their intersection is a singleton. These
{\em separating points} are the endpoints of unstable manifolds which
originate from saddle points. For later use we renumber the arcs $A_{ji}$
as $A_1,A_2,\ldots,A_s$ in counter-clockwise direction.

It is obvious that the number $s$ of separating points cannot be less
than the number of distinct zeros of $f$.
In order to get an upper bound of $s$ we assume that $f$ has $m$ distinct zeros
with multiplicities $\beta_1,\ldots,\beta_m$ and $k$ saddle points where
$f'$ has zeros of multiplicities $\alpha_1,\ldots,\alpha_k$, respectively.
Then we have
\[
\alpha_1 + \ldots + \alpha_k = m-1,\quad
\beta_1 + \ldots + \beta_m = n-1.
\]
The first equation follows from the well-known fact that the derivative
of a Blaschke product of order $n$ has exactly $n-1$ zeros in $\mathbb{D}$,
this time counting multiplicity.
From any saddle point $a_j$ exactly $\alpha_j+1$ rays emerge which belong
to the unstable manifold of $a_j$. Since any separating point must be
the endpoint of one such line, the total number $s$ of separating points
cannot be greater than $k+\beta_1 + \ldots + \beta_k = k+m-1$. Thus we
finally get
\[
m\le s \le m + k - 1.
\]
Examples show that both estimates are sharp.

It turns out that the global topological structure of the phase plot is
completely characterized by the sequence $S$ of integers, which associates
with any of the arcs $A_1,\ldots,A_s$ (in consecutive order) the number
of the corresponding zero. This sequence depends on the specific numbering
of the zeros and the arcs, but an appropriate normalization makes it unique.
For example, the Blaschke product depicted in Figure~24 is represented by
the sequence $S=(1,2,3,2,4,5,4,2)$.

Let us now return to Figure~1 again. Picturing once more that ``phase'' is
a substance emerging from sources at the zeros which can exit the domain
only at its boundary, is it then not quite natural that phase
plots of Blaschke products {\em must} look like they do?

And if you are asking yourself what ``natural'' means, then this is
already another question.

\section*{Concluding Remarks}

Phase plots result from splitting the information about the function
$f$ into two parts (phase and modulus), and one may ask why we do
not separate $f$ into its real and imaginary part.
One reason is that often zeros are of special interest; their presence
can easily be detected and characterized using the phase,
but there is no way to find these from the real or imaginary part
alone.

And what is the advantage of using $f/|f|$ instead of $\ln |f|$\,?
Of course, zeros and poles can be seen in the analytic landscape,
but they are much better represented in the phase plot. In fact there
is a subtle asymmetry between modulus and argument (respectively, phase).
For example, Theorem~\ref{t.rhp} has no counterpart for the modulus of
a function.%
\footnote{There is such a result for {\em outer functions}, but
it is impossible to see if a function is outer using only the boundary
values of its modulus.}

Since phase plots and standard domain coloring produce similar pictures,
it is worth mentioning that they are based on different concepts and have
a distinct mathematical background.

Recall that standard domain coloring methods use the complete values of
an analytic function, while phase plots depict only its phase. Taking into
account that phase can be considered as a periodization of the argument,
which is (locally) a harmonic function, reveals the philosophy behind phase
plots:
Analytic functions are considered as harmonic functions, endowed with a set
of singularities having a special structure.
Algebraically, phase plots forget about the linear structure of analytic
functions, while their multiplicative structure is preserved. 

This approach has at least two advantages. The first one is almost
trivial: phase has a small range, the unit circle, which allows one
visualizing all functions with one and the same color
scheme. Moreover, a one--dimensional color space admits a better resolution
of singularities. Mathematically more important is the existence of
a simple parametrization of {\em analytic} and {\em meromorphic} phase
plots by their boundary values and their singularities (Theorem~\ref{t.rhp}).
There is no such result for domain colorings of analytic functions.

The following potential fields of applications demonstrate that phase 
plots may be a useful tool for anyone working with complex--valued
functions. 

1. A trivial but useful application is {\em visual inspection} of
functions.
If, for example, it is not known which branch of a function is used in
a certain software, a glance at the phase plot may help.
In particular, if several functions are composed, software implementations
with different branch cuts can lead to completely different results.
You may try this with the {\sc Mathematica} functions ${\tt Log\,(Gamma)}$
and ${\tt LogGamma}$. Another useful exercise in teaching is to compare
the phase plots of $\exp(\log z)$ and $\log(\exp z)$.

2. A promising field of application is visual analysis and synthesis
of {\em transfer functions} in {\em systems theory} and {\em filter design}.
Since here the modulus (gain) is often more important than phase,
it is recommended to use the left color scheme of Figure~7.

3. Further potential applications lie in the area of {\em Laplace} and
complex {\em Fourier transforms}, in particular to the method of steepest
descent (or stationary phase).

4. Phase plots also allow to guess the {\em asymptotic behavior} of
functions (compare, for example, the phase plots of $\exp z$ and
$\sin z$), and to find {\em functional relations}. A truly challenging
task is to rediscover the functional equation of the Riemann Zeta--function 
from phase plots of $\zeta$ and $\Gamma$.

5. Complex {\em dynamical systems}, in the sense of iterated functions, have
been investigated by by Felix Huang \cite{Hua} and Martin Pergler \cite{Perg}
using domain coloring methods. The problem of scaling the modulus disappears
when using phase plots; see the pictures
of Fran\c{c}ois Labelle \cite{Lab} and Donald Marshall \cite{Mar}.

6. The utility of phase plots is not restricted to analytic functions.
Figure~25 visualizes the function
\[
h(z) := {\rm Im}\,\big({\mathrm{e}}^{-\frac{\mathrm{i}\pi}{4}}\,z^n\big)
+ \mathrm{i}\,{\rm Im}\,\big({\mathrm{e}}^{\frac{\mathrm{i}\pi}{4}}\,
(z-1)^n\big),
\]
with $n=4$. This is Wilmshurst's example \cite{Wil} of a harmonic
polynomial of degree $n$ having the maximal possible number of
$n^2$ zeros.
For background information we recommend the paper on
gravitational lenses by Dmitry Khavinson and Genevra Neumann
\cite{KhaNeu}.
\begin{center}\label{f.fig25}
\includegraphics[width=0.45\textwidth]{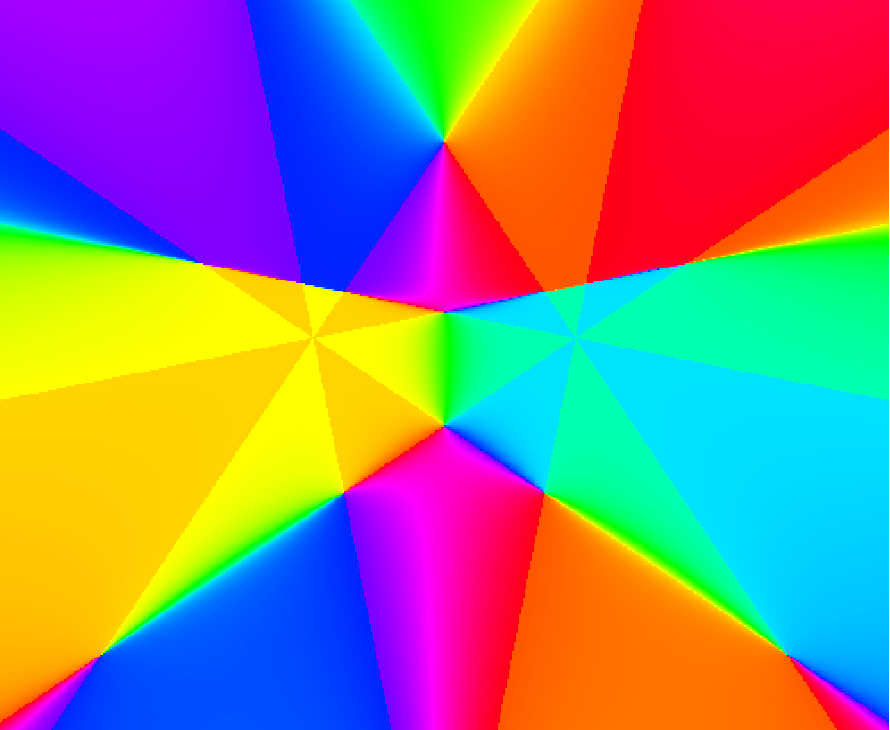}\\
Figure~25: A modified phase plot of Wilmshurst's example for $n=4$
\end{center}
To understand the construction of the depicted function it is important
to keep track of the zeros of its real and imaginary parts. In the figure
these (straight) lines are visualized using a modified color scheme
which has jumps at the points $1,\mathrm{i},-1$ and $-\mathrm{i}$ on the unit
circle.
\smallskip

Besides these and other concrete applications one important feature of phase
plots is their potential to bring up interesting questions and produce novel
ideas.
If you would like to try out phase plots on your own problems, you may start
with the following self-explaining {\sc Matlab} code:%
\footnote{This is a contribution to Nick Trefethen's project
of communicating ideas by exchanging ten-line computer codes.}
\smallskip

{\tt xmin=-0.5; xmax=0.5; ymin=-0.5; ymax=0.5; \\
xres = 400; yres = 400; \\
x = linspace(xmin,xmax,xres); \\
y = linspace(ymin,ymax,yres);} \\
$\mathbf{[}${\tt x,y}$\mathbf{]}$
{\tt = meshgrid(x,y); z = x+i*y; \\
f = exp(1./z); \\
p = surf(real(z),imag(z),0*f,angle(-f)); \\
set(p,'EdgeColor','none'); \\
caxis([-pi,pi]), colormap hsv(600) \\
view(0,90), axis equal, axis off}
\smallskip

Though the phase of a function is at least of the same importance as its
modulus, it has not yet been studied to the same extent as the latter. 
It is my conviction that {\em phase plots are problem factories}, which have
the potential to change this situation.
\smallskip

{\bf Technical Remark.}
All images of this article were created using {\sc Mathematica} and
{\sc Matlab}.
\smallskip

{\bf Acknowledgement.}
I would like to thank all people who supported me in writing this paper:
Gunter Semmler read several versions of the manuscript, his valuable ideas
and constructive criticism led to a number of significant improvements.
The project profited much from several discussions with Albrecht
B\"ottcher who also generously supported the production of a first
version. Richard Varga and George Csordas encouraged me during a
difficult period.
J\"{o}rn Steuding and Peter Meier kindly advised me in computing the
Riemann Zeta function.
Oliver Ernst tuned the language and eliminated some errors.

The rewriting of the paper would not have been possible without the
valuable comments and critical remarks of several referees.

Last but not least I would like to thank Steven Krantz and Marie
Taris for their kind, constructive, and professional collaboration
and an inspiring exchange of ideas.


\end{multicols}

Author's address: Elias Wegert,
Institute of Applied Mathematics,
Technical University Bergakademie Freiberg,
D-09596 Freiberg, Germany.
Email: {\tt wegert\symbol{64}math.tu-freiberg.de}

\end{document}